\documentclass[11pt,reqno]{article}

\usepackage{amscd, amssymb}
\usepackage[mathscr]{eucal}
\usepackage{mathrsfs}
\usepackage{amsfonts}
\usepackage{amsmath}
\usepackage{amsthm}
\usepackage{latexsym}

\theoremstyle{plain}
\newtheorem{theorem}[equation]{Theorem}

\newtheorem{lemma}[equation]{Lemma}

\theoremstyle{definition}
\newtheorem{definition}[equation]{Definition}

\newtheorem{example}[equation]{Example}
\newtheorem{remark}[equation]{Remark}

\makeatletter
\renewcommand{\subsection}{\@startsection{subsection}{2}{0pt}{-3ex
plus -1ex minus -0.2ex}{-2mm plus -0pt minus
-2pt}{\normalfont\bfseries}} \makeatother

\numberwithin{equation}{subsection}

\newcommand{\beq}{\begin{equation}\label}
\newcommand{\eeq}{\end{equation}}

\newcommand{\iso}{{\;\;\stackrel{_\sim}{\longrightarrow}\;\;}}

\newcommand{\vi}{${\sf {(i)}}\;$}
\newcommand{\vii}{${\sf {(ii)}}\;$}

\newcommand{\sset}{\subset}
\DeclareMathOperator{\Ker}{\mathtt{Ker}}
\DeclareMathOperator{\im}{\mathtt{Im}}

\def\map{\longrightarrow}

 \newcommand{\id}{{{\mathtt {Id}}}}

\newcommand{\into}{\,\,\hookrightarrow\,\,}

\newcommand{\Ad}{{\mathtt{{Ad}}^{\,}}}

\newcommand{\gr}{{\mathtt{{gr}^{\,}}}}

\newcommand{\lang}{\langle\langle\,}
\newcommand{\rang}{\,\rangle\rangle}
\newcommand{\lf}{\lfloor}
\newcommand{\rf}{\rfloor}

\newcommand{\rk}{{\mathtt{rk}}}

\newcommand{\Hom}{{\mathtt{Hom}}}

\newcommand{\la}{\lambda}
\newcommand{\om}{\omega}

\def\be{\beta}

\def\C{{\mathbb{C}}}
\def\R{{\mathsf{R}}}

\def\Z{{\mathbb{Z}}}

\def\B{{\mathsf B}}

\def\ZZ{{\mathsf{Z}}}

\def\A{{\mathsf{A}}}
\def\E{{\mathsf{E}}}
\newcommand{\g}{{\gamma}}
\newcommand{\G}{{\Gamma}}
\newcommand{\GG}{{\mathbf{\Gamma}_{n}}}
\newcommand{\hh}{{\mathsf{H}}}
\def\QQ{{\overline{Q}}}
\def\sy{{\mathcal{S}}}
\def\Up{{\Upsilon}}
\newcommand{\ve}{\varepsilon}
\newcommand{\ot}{\otimes}
\newcommand{\otot}{\otimes\cdots\otimes}

\begin{document}

\setlength{\parindent}{6mm}
\setlength{\parskip}{3pt plus 5pt minus 0pt}

\centerline{\Large {\bf Deformed preprojective algebras and symplectic}}
\vskip 2mm
\centerline{\Large {\bf reflection algebras for wreath products}}

\vskip 6mm
\centerline{\large {\sc Wee Liang Gan and Victor Ginzburg}}

\begin{abstract}
We determine the PBW deformations of the wreath product of
a symmetric group with a deformed preprojective algebra of an 
affine Dynkin quiver. In particular, we show that there is  
precisely one parameter which does not come from deformation
of the preprojective algebra. We prove that the PBW deformation 
is Morita equivalent to a corresponding symplectic reflection 
algebra for wreath product. 
\end{abstract}

\bigskip

\section{Introduction}
\setcounter{equation}{0}

\subsection{}
Deformed preprojective algebras were introduced by 
Crawley-Boevey and Holland in \cite{CBH}. We start by
recalling its definition.

Let $Q$ be a quiver, and denote by
$I$ the set of vertices of $Q$. 
The double $\QQ$ of $Q$ is the quiver obtained from $Q$ by  
adding a reverse edge $a^*: j\to i$ for each edge $a:i\to j$
in $Q$. For any edge  $a:i\to j$ 
in $\QQ$, we write its tail $t(a):=i$ and its head $h(a):=j$. 
Let $B:=\bigoplus_{i\in I} \C\,$,
and $E$ the vector space over $\C$ with 
basis formed by the set of edges $\{a\in \QQ\}$.
Thus, $E$ is a $B$-bimodule and $E=\bigoplus_{i,j\in I} E_{i,j}$,
where $E_{i,j}$ is spanned by the edges $a\in\QQ$ with $h(a)=i$
and $t(a)=j$. The path algebra of $\QQ$
is $\C\QQ := T_B E = \bigoplus_{n\geq 0} T^n_B E$, where
$T^n_B E = E\ot_B \cdots \ot_B E$ is the $n$-fold 
tensor product. 
The trivial path for the vertex $i$ is denoted by
$e_i$, an idempotent in $B$.
Let $r :=  \sum_{a\in Q}[a, a^*] \in T^2_B E$. 
For each $i\in I$, let $$r_i := e_ire_i=\sum_{\{a\in Q\,|\, h(a)=i\}} 
a\cdot a^*- \sum_{\{a\in Q\,|\, t(a)=i\}} a^*\cdot a \,.$$
For an element $\la\in B$,
we will write $\la= \sum_{i\in I}\la_i e_i$ where $\la_i\in\C$.

\begin{definition}
For each element $\la\in B$, the \emph{deformed preprojective 
algebra} of $Q$ is the quotient algebra
$$ \Pi_\la := \frac{\C\QQ}{\lang r-\la \rang}
=  \frac{\C\QQ}{\lang r_i-\la_i e_i \rang_{i\in I}}\,\,,$$
where $\lang\ldots\rang$ is the two-sided ideal generated by
the indicated elements.
\end{definition}

The algebra $\Pi_0$ is called the preprojective algebra of $Q$.
Note that the grading on $\C\QQ$ induces a filtration on
$\Pi_{\la}$ and there is a natural map $\Pi_0 \map \gr\Pi_{\la}$.
The ``PBW'' theorem for $\Pi_{\la}$ proved in \cite[Cor. 3.6]{CBH}
says that this map is an isomorphism when $Q$ is
an affine Dynkin quiver (of type ADE).

\subsection{}
The main construction of this paper is a one-parameter deformation 
of the wreath product $\Pi_{\la}^{\ot n} \# S_n$, where 
$n$ is an integer greater than $1$, the superscript
$\ot n$ means $n$-fold tensor product over $\C$, 
and $S_n$ is the symmetric group on $n$ objects.

To state the definition, we will use the following notations. 
Let $n$ be a positive integer. The element $s_{ij}\in S_n$ is
the transposition $i\leftrightarrow j$.
Let $\B := B^{\ot n}$. For any $\ell\in [1,n]$, define the 
$\B$-bimodules
$$ \E_{\ell}:= B^{\ot (\ell-1)}\ot E\ot B^{\ot (n-\ell)}
\qquad \mathrm{and}\qquad \E := 
\bigoplus_{1\leq\ell\leq n} \E_\ell\,.$$
The natural inclusions $\E_\ell \into 
B^{\ot (\ell-1)}\ot T_B E\ot B^{\ot (n-\ell)} \sset (T_B E)^{\ot n}$ 
induce canonical identifications
$T_\B\E_\ell = B^{\ot (\ell-1)}\ot T_B E\ot B^{\ot (n-\ell)}$ and
a surjective morphism $\Up:T_{\B}\E \map (T_B E)^{\ot n}$.
Given two elements $\ve\in\E_\ell$ and $\ve'\in\E_m$ of the form
\begin{equation} \label{eqnve}
\ve= e_{i_1}\ot e_{i_2}\otot a\otot h(b) \otot e_{i_n}\,,
\end{equation}
\begin{equation} \label{eqnve'}
\ve'= e_{i_1}\ot e_{i_2}\otot t(a)\otot b\otot e_{i_n}\,,
\end{equation}
where $\ell\neq m$, $a,b\in \QQ$ and $i_1, \ldots, i_n \in I$,
we define their ``commutator''
\begin{eqnarray}
\lf\ve,\ve'\rf &:=& 
(e_{i_1}\otot a\otot h(b) \otot e_{i_n})
(e_{i_1}\otot t(a)\otot b \otot e_{i_n}) \nonumber\\ & &
-(e_{i_1}\otot h(a)\otot b \otot e_{i_n})
(e_{i_1}\otot a\otot t(b) \otot e_{i_n}) . \nonumber
\end{eqnarray}
Note that $\lf\ve,\ve'\rf$ is an element in $T^2_\B \E$.
The kernel of $\Up$ is the two-sided ideal
generated by all elements of the form $\lf\ve,\ve'\rf$.

\begin{definition} \label{dea}
Let $n>1$ be an integer.
For any $\la\in B$ and $\nu\in \C$, define the algebra
$\A_{n,\la,\nu}$ to be the quotient of $T_{\B}\E\,\# S_n$
by the following relations.
\begin{itemize}
\item[\vi] 
For any $i_1,\ldots, i_n\in I$ and $\ell\in [1,n]$:
$$e_{i_1}\otot (r_{i_\ell}-\la_{i_\ell}e_{i_\ell}) 
\otot e_{i_n} = \nu \sum_{\{ j\neq\ell \,|\, i_j=i_\ell\}}
(e_{i_1}\otot e_{i_\ell} \otot e_{i_n})s_{j\ell}\,.$$
\item[\vii]
For any $\ve, \ve'$ of the form (\ref{eqnve})--(\ref{eqnve'}):
\[  \lf\ve,\ve'\rf = \left\{ \begin{array}{ll}
\nu(e_{i_1}\otot h(a)\otot t(a)\otot e_{i_n})s_{\ell m}
& \textrm{if $a = b^*$ \& $b\in Q$}\,,\\
- \nu(e_{i_1}\otot h(a)\otot t(a)\otot e_{i_n})s_{\ell m}
& \textrm{if $b = a^*$ \& $a\in Q$}\,,\\
0 & \textrm{else}\,.  \end{array} \right.  \]
\end{itemize}
\end{definition}

It is easy to see that $\A_{n,\la,\nu}$ does not depend on
the orientation of $Q$, cf. \cite[Lemma 2.2]{CBH}.
Moreover, $\A_{n,\la,0} = \Pi^{\ot n}_{\la} \# S_n$.
The grading on $T_{\B}\E$ induces a filtration
on $\A_{n,\la,\nu}$ and there is a natural map
$ \Pi^{\ot n}_{0}\# S_n=\A_{n,0,0} \map\gr\A_{n,\la,\nu}$.
We will prove in \S \ref{s2} that, when $Q$ is affine Dynkin
of type ADE, this map is an isomorphism, and
any ``PBW deformation'' of $\Pi^{\ot n}_{0}\# S_n$ must be 
of the form $\A_{n,\la,\nu}$.

\begin{remark}
Regarding $\nu$ as a formal variable,
we obtain a one-parameter \emph{formal} deformation of 
$\Pi^{\ot n}_{\la} \# S_n$; see \cite{EO}.
\end{remark}

\subsection{} \label{s13}
The motivation to study $\A_{n,\la,\nu}$ comes from \cite{EG},
in which Etingof and Ginzburg introduced the so-called symplectic
reflection algebras for wreath products, cf. also \cite{GS}.
We recall its definition.

Let $L$ be a $2$-dimensional vector space equipped with a 
nondegenerate symplectic form $\om_L$.
Let $V := L^{\oplus n}$ and $\om : =\om_{L}^{\oplus n}$.
Let $\G$ be a finite subgroup of $Sp(L)$ and
$\GG :=S_{n} \ltimes \G^{n} \subset Sp(V)$.
Denote by $\ZZ\G$ the center of the group algebra $\C[\G]$.
Given $\g\in \G$, write $\g_{i} \in \GG$ for $\g$ placed in the
$i$-th factor $\G$.
An element $s\in \GG$ is called a symplectic reflection if
$\rk(\id-s)=2$. According to \cite[(11.1)]{EG},
there are two types of symplectic reflections in $\GG$:

\medskip
(S) The elements $s_{ij}\g_{i}\g_{j}^{-1}$, where $i,j \in [1,n]$
and $\g \in \G$.

($\G$) The elements $\g_{i}$, where $i \in [1,n]$ and 
$\g \in \G\smallsetminus \{1\}$.
\medskip

The group $\GG$ acts on the set $\sy$ of symplectic reflections
by conjugation. The set of elements of type (S) form a single 
$\GG$-conjugacy class, while the elements of type ($\G$) 
for $\g$ in each $\G$-conjugacy class form a $\GG$-conjugacy class.
Thus, we may identify an $\Ad\GG$-invariant function
$c:\sy\map\C\,:\,s\mapsto c_s$ with an element 
$k\cdot 1 + \sum_{\g \in \G\smallsetminus \{1\}} 
c'_{\g}\cdot\g \in \ZZ\G$, where $k$ is the value of $c$ on
elements of type (S) and $c'_{\g}$ is the value of $c$ on
the elements $\g_i$ of type ($\G$).
For each $s\in\sy$, write $\om_s$ for the bilinear form 
on $V$ which coincides with $\om$ on $\im(\id-s)$
and has $\Ker(\id-s)$ as its radical.

\begin{definition}
For any $t\in\C$ and $c\in \ZZ\G$, the \emph{symplectic reflection
algebra} $\hh_{t,c}(\GG)$ is defined to be the quotient algebra 
$(TV\#\GG)/\lang [u,v]-\kappa(u,v) \rang_{u,v \in V} $, where
$$\kappa\,:\,V\ot V\map \C[\GG]\,:\, (u,v)\mapsto
 t\cdot\om(u,v)\cdot 1 +\sum_{s\in \sy}
c_s\cdot\om_s (u,v)\cdot s  \,.$$
\end{definition}

We will construct in \S \ref{s3} a Morita equivalence between 
$\hh_{t,c}(\GG)$ and $\A_{n,\la,\nu}$, where the quiver $Q$ is
associated to $\G$ via the \emph{McKay correspondence}.

\begin{remark} \label{rmk10}
When $\G=\{1\}$, the algebra $\hh_{t,c}(\GG)=\hh_{t,k}(S_n)$ is the
rational Cherednik algebra of type $A_{n-1}$, cf.
\cite{EG} and Lemma \ref{leh} below. 
In this case, the quiver $Q$ is the affine Dynkin
quiver of type $A_0$, and we have $\A_{n,\la,\nu}=\hh_{t,k}(S_n)$
where the parameters $\la=t$ and $\nu=k/2$.
\end{remark}

\section{PBW deformation} \label{s2}

\subsection{}
We will define here what we mean by PBW deformations of 
$\Pi^{\ot n}_0\# S_n$.

Denote by $\lf \E_\ell, \E_m\rf$ the sub-$\B$-bimodule of 
$T^2_{\B} \E$ spanned by all elements of the form $\lf\ve, \ve'\rf$
with $\ve\in\E_\ell$, $\ve'\in\E_m$.
Let $R$ be the sub-$B$-bimodule of $T^2_B E$ spanned by $r_i$
for $i\in I$, and let
$$ \R_{\ell}:= B^{\ot (\ell-1)}\ot R \ot B^{\ot (n-\ell)}
\,,\quad \R:= \big( \bigoplus_{1\leq\ell\leq n}\R_\ell \big)\oplus
\big(\bigoplus_{1\leq\ell <m\leq n}\lf \E_\ell,\E_m\rf \big)
\sset T^2_\B \E\,.$$
Furthermore, let $K:=\B\# S_n$ and $M:= \E\ot \C[S_n]$.
Note that $M$ is a $K$-bimodule, where the left action of $S_n$
on $M$ is the diagonal one. We have: $T_\B \E\,\# S_n = T_K M$.
Let $U:=\R\ot \C[S_n] \sset T^2_K M$. 
For any $K$-bilinear map $\be:U \map K$, define the algebra
$$\A_\be := \frac{T_\B \E\,\# S_n}{\lang x-\be(x) \rang_{x\in U}} 
= \frac{T_K M}{\lang x-\be(x) \rang_{x\in U}}\,. $$
The grading on $T_\B \E\,\# S_n$ induces a filtration on
$\A_\be$ and there is a natural map $\Pi^{\ot n}_0\# S_n
= \A_0 \map \gr\A_\be$.

\begin{definition}
The algebra $\A_\be$ is a \emph{PBW deformation} of $\Pi^{\ot n}_0
\# S_n$ if $\gr\A_\be = \Pi^{\ot n}_0\# S_n$.
\end{definition}

\begin{remark}
When $n=1$ and $Q$ is affine Dynkin of type ADE, 
any $B$-bilinear map $\be: R \map B$ gives a PBW 
deformation by \cite[Cor. 3.6]{CBH}.
\end{remark}

\subsection{}
The first main result of this paper is the following.

\begin{theorem} \label{Thm1}
Let $n>1$ and assume $Q$ is an affine Dynkin quiver of type ADE. 
The algebra $\A_\be$ is a PBW deformation of 
$\Pi^{\ot n}_0 \# S_n$ if and only if $\A_\be = \A_{n,\la,\nu}$
for some $\la\in B$, $\nu\in \C$.
\end{theorem}

\begin{proof}
When $Q$ is of type $A_0$, the Theorem follows from 
\cite[Theorem 1.3]{EG}; see Remark \ref{rmk10}. 
Hence, we may assume that $Q$ has no edge-loop.
Given any $x\in U$, write $\be(x)$ in the form   
$$\be(x)=\sum \be_{\sigma}^{j_1\cdots j_n} (x)
e_{j_1}\otot e_{j_n} \cdot\sigma $$ 
where $\be_{\sigma}^{j_1\cdots j_n} (x) \in \C$ and
the sum is taken over all vertices $j_1, \ldots, j_n \in I$ 
and permutations $\sigma\in S_n$.

First, we find the constraints on 
$\be_{\sigma}^{j_1\cdots j_n}$ so that $\be$ is
$K$-bilinear.
The right $S_n$-linearity of $\be$ means that
the $\be_{\sigma}^{j_1\cdots j_n}$'s are determined by how
they are defined on $\R$.
The left $S_n$-linearity of $\be$ means that 
for any $\sum k_1\ot k_2 \otot k_n \in \R$ and $\tau\in S_n$,
we have
$$\be_{\sigma}^{j_1\cdots j_n}(\sum k_1\ot k_2\otot k_n)   
= \be_{\tau\sigma\tau^{-1}}^{j_{\tau(1)}\cdots j_{\tau(n)}}
(\sum k_{\tau(1)}\ot k_{\tau(2)}\otot k_{\tau(n)})\,.$$

Consider any element $k_1 \otot k_n \in \R_\ell$
with $k_1, k_2, \ldots \in I$ except for $k_\ell = r_i \in R$.
By the $\B$-bilinearity of $\be$, we must have
$\be_{\sigma}^{j_1\cdots j_n}(k_1 \otot k_n) = 0$ if
$j_1 \neq k_1$, or $j_2 \neq k_2, \ldots$, or $j_\ell \neq i,\ldots$ 
or $j_n \neq k_n$, or $k_{\sigma(p)}\neq k_p$ for any $p$.   

Consider any element $\lf\ve,\ve'\rf$ with $\ve,\ve'$
of the form (\ref{eqnve})--(\ref{eqnve'}).
By the $\B$-bilinearity of $\be$, we must have
$\be_{\sigma}^{j_1\cdots j_n}(\lf\ve,\ve'\rf) = 0$ if
$j_1 \neq i_1$, or $j_2 \neq i_2, \ldots$, or $j_\ell \neq h(a)$
or $j_m\neq h(b)$, or $i_{\sigma(p)}\neq i_p$ for any $p$, or
($\ell,m$ in different cycles of $\sigma$), or
($\ell,m$ in same cycle of $\sigma$ and $t(b)\neq h(a)$ or
$t(a)\neq h(b)$).
Note, in particular, that $\be(\lf\ve,\ve'\rf)=0$ 
if $t(b)\neq h(a)$ or $t(a)\neq h(b)$.

Next, it is known that the $B$-algebra $\Pi_0$ is Koszul; 
this was proved in \cite[Theorem 7.2]{Gr} when $Q$ is affine Dynkin 
of type $A$ and in \cite[Theorem 1.9]{MV} when $Q$ is affine Dynkin 
of type $D$ or $E$.
Thus, $\Pi^{\ot n}_0 \# S_n$ is a Koszul $K$-algebra, and
so by \cite[Lemma 3.3]{BG}, $\A_\be$ is a PBW deformation if and only if
$\be\ot\id=\id\ot\be$ on $\big( \R\ot_\B \E\big) \cap 
\big( \E\ot_\B \R\big)$. Here, the equality
takes place in $M$ while the intersection takes place in $T^3_K M$.

Observe that $\big( \R\ot_\B \E\big) \cap\big( \E\ot_\B \R\big)$ 
is spanned by the following two types of elements:

(1) For any $\ve\in\E_\ell$, $\eta\in\E_m$ and $\zeta\in\E_r$ 
of the form
\begin{equation}\label{eq5}
 \ve = e_{i_1}\otot a \otot h(b) \otot h(c) \otot e_{i_n}\,,
\end{equation}
\begin{equation}\label{eq6}
 \eta= e_{i_1} \otot t(a)\otot b \otot h(c) \otot e_{i_n}\,,
\end{equation}
\begin{equation}\label{eq7}
 \zeta= e_{i_1}\otot t(a)\otot t(b) \otot c \otot e_{i_n}\,,
\end{equation}
where $a,b,c\in \QQ$ and $i_1, \ldots, i_n\in I$, we have the
element:
$$ \lf\ve,\eta\rf\zeta -\lf\ve,\zeta\rf\eta +\lf\eta,\zeta\rf\ve
= \ve\lf\eta,\zeta\rf -\eta\lf\ve,\zeta\rf +\zeta\lf\ve,\eta\rf \,.$$
Here, in the second term of the left hand side, the $\eta$ is
actually the $\eta$ of (\ref{eq6}) whose $r$-th entry is $t(c)$
instead of $h(c)$, and the $\zeta$ is actually the $\zeta$ of
(\ref{eq7}) whose $m$-th entry is $h(b)$ instead of $t(b)$. 
Throughout, we shall use similar convention.

(2) We use a similar convention as above. 
For any $x=\sum \ve_i \eta_i\in \R_\ell$ (where 
$\ve_i, \eta_i\in \E_\ell$) and $\zeta\in \E_m$, we have the element:
$$\big(\sum \ve_i \eta_i\big)\zeta-\sum\lf \ve_i,\zeta\rf\eta_i 
= \sum \ve_i\lf\eta_i,\zeta\rf+\zeta\big(\sum \ve_i \eta_i\big)\,.$$

We now see what the equation $\be\ot \id = \id \ot\be$ says
when applied to elements of type (1). We have:
$$ \be(\lf\ve,\eta\rf)\zeta -\be(\lf\ve,\zeta\rf)\eta 
+\be(\lf\eta,\zeta\rf)\ve = \ve\be(\lf\eta,\zeta\rf) 
-\eta\be(\lf\ve,\zeta\rf) +\zeta\be(\lf\ve,\eta\rf) \,.$$
Recall that this is an equality of elements in $M$.
The edges which appear in $\be(\lf\ve,\eta\rf)\zeta$,
$\be(\lf\ve,\zeta\rf)\eta$, and $\be(\lf\eta,\zeta\rf)\ve$
are, respectively, $c$, $b$, and $a$; the same for
$\zeta\be(\lf\ve,\eta\rf)$, $\eta\be(\lf\ve,\zeta\rf)$,
and $\ve\be(\lf\eta,\zeta\rf)$.
It follows that for equality, we must have 
$ \be(\lf\ve,\eta\rf)\zeta = \zeta\be(\lf\ve,\eta\rf)$ 
when $c\neq a,b$. Hence, we deduce that
$\be_{\sigma}^{j_1\cdots j_n}(\lf\ve,\eta\rf) = 0$ 
if $\sigma\neq s_{\ell m}$.
Moreover, we also have $$\be_{\sigma}^{j_1\cdots j_r \cdots j_n}
(\lf \ve,\eta\rf)
= \be_{\sigma}^{j_1\cdots j'_r \cdots j_n}(\lf\ve,\eta\rf)$$ 
where $r\neq \ell,m$ and $j_r= h(c)$, $j'_r=t(c)$ (for any edge
$c\in\QQ$).
Here, by our convention, the $r$-th entry of
$\lf \ve,\eta\rf$ on the left hand side is $h(c)$ 
while the $r$-th entry of $\lf \ve,\eta\rf$ on the right 
hand side is $t(c)$.

Now we see what the equation $\be\ot \id = \id \ot\be$ says
when applied to elements of type (2). We get:
\begin{equation} \label{eq8}
\be\big(\sum \ve_i \eta_i\big)\zeta -
\sum \be\big(\lf \ve_i,\zeta\rf \big)\eta_i
= \sum \ve_i\be\big(\lf\eta_i,\zeta\rf \big) +
\zeta\be\big(\sum \ve_i \eta_i\big)\,.
\end{equation}
Hence,
$$\be_{1}^{j_1\cdots j_m \cdots j_n}(x)
= \be_{1}^{j_1\cdots j'_m \cdots j_n}(x)$$ 
where $m\neq\ell$ and $j_m= h(c)$, $j'_m=t(c)$ (for any edge $c\in\QQ$). 
Here, by our convention, the $m$-th entry of $x$  on the left hand 
side is $h(c)$ while the $m$-th entry of $x$ on the right
hand side is $t(c)$.
Moreover, we have $\be_{\sigma}^{j_1\cdots j_n}(x) =0$ if
$\sigma$ is not $1$ or $s_{\ell m}$ for any $m$.

Let $\sigma=s_{\ell m}$. Taking $\zeta$ to be $\ve'$ of
(\ref{eqnve'}) and $x$ to be $e_{i_1}\otot r_{i_\ell}\otot e_{i_n}$,
we deduce from (\ref{eq8}) that
$\be_\sigma^{i_1\cdots i_n}(\lf\ve,\ve'\rf) 
  = \be_\sigma^{i_1\cdots i_n} (x)$ if $b\in Q$ and $a=b^*$;
$\be_\sigma^{i_1\cdots i_n}(\lf\ve,\ve'\rf) 
  = -\be_\sigma^{i_1\cdots i_n}(x)$ if $a\in Q$ 
and $b=a^*$;
and $\be_\sigma^{i_1\cdots i_n}(\lf\ve,\ve'\rf) = 0$ otherwise.
\end{proof}

\begin{remark}
Let $Q$ be a connected quiver.
It is known from \cite{Gr} and \cite{MV} that if $Q$ is not
a finite Dynkin quiver, then the preprojective 
algebra $\Pi_0$ is Koszul. In this case, note that
$(R\ot_B E) \cap(E\ot_B R) =0$ when $Q$ has more than one edge,
and so by \cite{BG}, the deformed 
preprojective algebra $\Pi_\la$ is PBW for any $\la\in B$;
moreover, assuming furthermore that $Q$ has no edge-loop,
Theorem \ref{Thm1} is still true by same proof as above.
\end{remark}

\subsection{}
We end this section with some comments. First, $T_\B \E$ is the
path algebra of the product quiver 
$\QQ \times \cdots \times \QQ$
whose vertex set is $I\times\cdots\times I$ and edge set is 
$\bigcup I\times\cdots \times \QQ\times\cdots\times I$.

Next, let us consider the relations \vii in Definition \ref{dea}
for $n=2$.

\begin{example}\label{ex2}
When $n=2$, the relations \vii in Definition \ref{dea}
means that, for any edge $a\in Q$:
$$ (a^*\ot h(a))(h(a)\ot a)-(t(a)\ot a)(a^*\ot t(a))
= \nu \cdot(t(a)\ot h(a))s_{12}\,; $$
and for any edges $a,b\in\QQ$ with $a\neq b^*$ or $b\neq a^*$: 
$$ (a\ot h(b))(t(a)\ot b)-(h(a)\ot b)(a\ot t(b)) =0\,.$$
\end{example}

\section{Morita equivalence} \label{s3}

\subsection{}
We will give in the following lemma a more explicit presentation
of the algebra $\hh_{t,c}(\GG)$ by generators and relations.
Given any element $u\in L$ and $i\in [1,n]$, 
we will write $u_i\in V$ for $u$ placed in the $i$-th factor $L$.
Recall from \S \ref{s13} that $c = k\cdot 1 + \sum_{\g \in \G
\smallsetminus \{1\}}c'_{\g}\cdot\g \in \ZZ\G$.
From now on, we fix a basis $\{x, y\}$ of $L$ with $\om_L(x,y)=1$.

\begin{lemma} \label{leh}
The algebra $\hh_{t,c}(\GG)$ is the quotient
of $TV\#\GG$ by the following relations:
\begin{itemize}
\item[{\emph{(R1)}}] 
For any $i\in [1,n]$: $$[x_{i}, y_{i}]
= t\cdot 1+ \frac{k}{2} \sum_{j\neq i}\sum_{\g\in\G}
s_{ij}\g_{i}\g_{j}^{-1} + \sum_{\g\in\G\smallsetminus\{1\}}
c'_{\g}\g_{i}\,.$$
\item[{\emph{(R2)}}] 
For any $u,v\in L$ and $i\neq j$:
$$[u_{i},v_{j}]= -\frac{k}{2} \sum_{\g\in\G} \omega_{L}(\g u,v)
s_{ij}\g_{i}\g_{j}^{-1} \,.$$
\end{itemize}
\end{lemma}

\begin{proof}
We first consider symplectic reflections of type (S). 
Let $s= s_{ij}\g_{i}\g_{j}^{-1}$.
If $u \in V$, then
$(u-su)/2 \in \im(\id-s)$ and $(u+su)/2 \in \Ker(\id-s)$. Thus,
for any $u,v \in V$, we have
$\omega_{s}(u,v) = \omega(u-su, v-sv)/4
= \omega(u,v)/2 - \omega(u,sv)/2$.
In particular,
\[ \begin{array} {rcll}
\omega_{s}(x_{i}, y_{i}) &=& 1/2 \,,  &  \nonumber \\
\omega_{s}(x_{l}, y_{l}) &=& 0 \quad & 
 \textrm{for $l \neq i,j$} \,,\nonumber\\
\omega_{s}(u_{i}, v_{j}) &=& -\omega_{L}(u,\g^{-1}v)/2 \quad & 
\textrm{for any $u,v \in L$ } \,,\nonumber\\
\omega_{s}(u_{l}, v_{m}) &=& 0 \quad &  \textrm{for any 
$u,v \in L$, $l\neq m, i, j$} \,.\nonumber
\end{array} \]

We next consider symplectic reflections of type ($\G$).
Let $s=\g_{i}$. Then $\im(\id-s)= \{u_{i}|u\in L\}$, and
$\Ker(\id-s)$ is spanned by $u_{j}$ where $u\in L$, $j \neq i$.
Thus,
\[ \begin{array}{rcll}
\omega_{s}(x_{i}, y_{i})&=& 1 \,, & \nonumber\\
\omega_{s}(x_{l}, y_{l})&=& 0 \quad &
\textrm{for $l \neq i$} \,, \nonumber\\
\omega_{s}(u_{l}, v_{m})&=& 0 \quad & \textrm{for any $u,v \in L$, 
$l \neq m, i$ } \,. \nonumber
\end{array} \]
\end{proof}

\subsection{}
Let us recall the classical McKay correspondence. Given the finite 
subgroup $\G\sset Sp(L)$, we shall write its irreducible 
representations as $N_i$, $i\in I$. Consider the quiver with
vertex set $I$ and whose number of edges from $i$ to $j$ is 
the multiplicity of $N_i$ in $L\ot N_j$. 
This quiver is the double of an affine Dynkin 
quiver $Q$ of type ADE. This construction 
gives a bijection between conjugacy classes of finite subgroups
of $Sp(L)$ and affine Dynkin diagrams of type ADE.

\subsection{}
Following \cite[\S 4]{CB} and \cite[\S 3]{CBH}, we now
define idempotent elements $f_i$ and $f$ in the group algebra $\C\G$.
For each $i\in I$, let $\delta_{i}$ be the dimension 
of the irreducible representation $N_{i}$. We fix an isomorphism
$\mathbb{C}\G \simeq \bigoplus_{i\in I} 
\mathrm{Mat}(\delta_{i}\times\delta_{i})$.
Let $E^{i}_{p,q}$ ($1\leq p,q \leq \delta_{i}$)
be the element of $\mathbb{C}\G$ with
$1$ in the $(p,q)$-entry of the matrix for the $i$-th summand 
and zero elsewhere.
Let $f_{i}$ be the idempotent $E^{i}_{1,1}$, and let
$f = \sum_{i\in I} f_{i}$. 

Note that, in the algebra $\mathbb{C}[\G^{n}]=
(\mathbb{C}\G)^{\otimes n}$, we have $$f^{\ot n}
= \sum_{i_1, \ldots, i_n \in I} f_{i_1}\otot f_{i_n}$$ and
\begin{align} \label{eqn11}
 \sum_{i_{1}, p_{1}, \ldots, i_{n}, p_{n}} &
(E^{i_{1}}_{p_{1},1} \otimes \cdots \otimes E^{i_{n}}_{p_{n},1})
f^{\otimes n} (E^{i_{1}}_{1,p_{1}} \otimes \cdots \otimes
E^{i_{n}}_{1,p_{n}})  \nonumber\\
= \sum_{i_{1}, p_{1}, \ldots, i_{n}, p_{n}} &
E^{i_{1}}_{p_{1},1}E^{i_{1}}_{1,p_{1}} \otimes \cdots \otimes
E^{i_{n}}_{p_{n},1}E^{i_{n}}_{1,p_{n}}
\quad =\, 1^{\otimes n} \,.
\end{align}

\subsection{} \label{ss34}
We state here some observations which we will use later.
First, we have an isomorphism
\begin{equation} \label{eqn14}
\B \iso f^{\ot n}\C[\G^n]f^{\ot n} = \bigoplus_{i_1,\ldots,
i_n\in I} \C\cdot f_{i_1}\otot f_{i_n}
\end{equation}
defined by
$$ e_{i_1}\otot e_{i_n} \mapsto f_{i_1}\otot  f_{i_n} \,.$$
Now $V \ot \C[\G^{n}]$ is a $\C[\G^n]$-bimodule, where the 
left action is the diagonal one. We have:
\begin{align} \label{eqn15}
& f^{\ot n} (V \ot \mathbb{C}[\G^{n}]) f^{\ot n} \nonumber\\
= & \bigoplus_{i_{1}, \ldots, j_{n}}
(f_{i_{1}} \otimes \cdots \otimes f_{i_{n}})
(L^{\oplus n} \otimes \underbrace{
\mathbb{C}\G \otimes \cdots \otimes \mathbb{C}\G}_n )
(f_{j_{1}} \otimes \cdots \otimes f_{j_{n}}) \nonumber\\
=\, & \bigoplus_{l=1}^{n} \bigoplus_{i_{1}, \ldots, j_{n}}
\Hom_{\G}(N_{i_{1}}, N_{j_{1}}) \otot
\Hom_{\G}(N_{i_{l}}, L\otimes N_{j_{l}})
\otimes \cdots \otimes \Hom_{\G}(N_{i_{n}}, N_{j_{n}}) \nonumber\\
\simeq\, & \E \,.
\end{align}
It follows from (\ref{eqn14})--(\ref{eqn15}) that
\begin{equation} \label{eqn16}
 f^{\otimes n} T_{\mathbb{C}[\G^{n}]} (V \ot \mathbb{C}[\G^{n}])
f^{\otimes n} \simeq T_\B \E
\end{equation}
and
\begin{equation} \label{eqn17}
 f^{\otimes n} (TV\#\GG) f^{\otimes n}
= f^{\otimes n} (T_{\mathbb{C}[\G^{n}]} (V \otimes
\mathbb{C}[\G^{n}]) \# S_{n}) f^{\otimes n}
\simeq T_\B \E \, \# S_{n} \,.
\end{equation}

\subsection{}
By (\ref{eqn11}), the algebra $\hh_{t,c}(\GG)$ is Morita
equivalent to the algebra $f^{\ot n}\hh_{t,c}(\GG)f^{\ot n}$.
By (\ref{eqn17}), $f^{\ot n}\hh_{t,c}(\GG)f^{\ot n}$ is  
isomorphic to a quotient of $T_\B \E \, \# S_{n}$.

\begin{remark}
When $n=1$, there is no parameter $k$ and $\hh_{t,c}(\GG)$ is 
the algebra denoted by $\mathscr{S}^\la$ in \cite{CBH},
with $\la_i$ being the trace of $t\cdot 1+\sum_{\g\neq 1}c'_\g \g$
on $N_i$; cf. Lemma \ref{leh}.
In this case, it was proved in \cite[Theorem 3.4]{CBH} 
that $f\mathscr{S}^\la f \simeq \Pi_\la$.
\end{remark}

The second main result of this paper is the following
generalization of \cite[Theorem 3.4]{CBH}; cf. also 
\cite[Theorem 4.12]{CB}.

\begin{theorem}
When $n>1$, there is an isomorphism 
$f^{\ot n}\hh_{t,c}(\GG)f^{\ot n} \simeq \A_{n,\la,\nu}$,
where $\la_i$ is the trace of $t\cdot 1+\sum_{\g\neq 1}c'_\g \g$
on $N_i$ and $\nu= \frac{k|\G|}{2}$.
\end{theorem}

\begin{proof}
When $\G = \{1\}$, the Theorem is trivial by Remark \ref{rmk10}.
Hence, we may assume $\G \neq \{1\}$, and so $Q$ is not of type
$A_0$. 

Denote by $\zeta: \mathbb{C} \rightarrow L\otimes L$ the linear
map that sends $1$ to $y \ot x - x \ot y$.
By Lemma 3.2 of \cite{CBH} and its proof,  
for each edge $a \in Q$, there are $\G$-equivariant monomorphisms
$$\theta_{a}: N_{t(a)} \map L\ot N_{h(a)} \qquad \textrm{and}\qquad 
\phi_{a}: N_{h(a)} \map L\ot N_{t(a)} $$
such that for each vertex $i$, we have
$$ \sum_{a\in Q, h(a)=i}(\id_L \otimes \theta_{a})\phi_{a} -
\sum_{a\in Q, t(a)=i}(\id_L \otimes \phi_{a})\theta_{a} =
-\delta_{i}(\zeta \otimes \id_{N_i})$$
as maps from $N_i$ to $L\ot L\ot N_i$, and such that 
$$(\omega_{L}\ot \id_{N_{t(a)}})(\id_L \ot \phi_{a})\theta_{a} =
-\delta_{h(a)}\id_{N_{t(a)}}$$ and
$$(\omega_{L}\ot \id_{N_{h(a)}})(\id_L \ot \theta_{a})\phi_{a} =
\delta_{t(a)}\id_{N_{h(a)}}\,.$$
Moreover, the $\theta_a, \phi_a$ ($a\in Q$) combine to give
a basis for each of the spaces $\Hom_\G (N_i, L\ot N_j)$. 
By (\ref{eqn14})--(\ref{eqn17}), we have an isomorphism 
$T_\B \E \#S_n \iso f^{\ot n}(TV\#\,\GG)f^{\ot n}$ 
such that
$$ e_{i_1}\ot e_{i_2} \otot e_{i_n} \cdot \sigma
\mapsto f_{i_1}\ot f_{i_2} \otot f_{i_n}\cdot \sigma\,,$$
$$ e_{i_1}\ot e_{i_2} \otot a \otot e_{i_n} \cdot \sigma \mapsto  
f_{i_1}\ot f_{i_2} \otot \phi_a \otot f_{i_n}\cdot \sigma\,,$$
$$ e_{i_1}\ot e_{i_2} \otot a^* \otot e_{i_n} \cdot \sigma\mapsto 
f_{i_1}\ot f_{i_2} \otot \theta_a \otot f_{i_n}\cdot \sigma\,,$$
for any $i_1, \ldots, i_n \in I$, $a \in Q$, and $\sigma \in S_n$.

Denote by $J$ the subspace of $T V\,\#\GG$ spanned by elements
of the form $[u,v] -\kappa(u,v)$ with $u,v \in V$. The algebra
$\hh_{t,c}(\GG)$ is the quotient of $T V\,\#\GG$ by the two sided
ideal generated by $J$. Thus, $f^{\ot n} \hh_{t,c}(\GG) f^{\ot n}$
is the quotient of $f^{\ot n}(TV\,\#\GG)f^{\ot n}$ by the ideal
$ f^{\ot n}(TV\#\,\GG) J (TV\#\,\GG) f^{\ot n}$.
By (\ref{eqn11}), we have 
$$ f^{\ot n}(TV\#\,\GG) J (TV\#\,\GG) f^{\ot n}
=  f^{\ot n}(TV\#\,\GG) f^{\ot n} \C[\G^n]
J \C[\G^n]f^{\ot n}(TV\#\,\GG) f^{\ot n} \,.$$
Hence, $f^{\ot n} \hh_{t,c}(\GG) f^{\ot n}$ is the quotient 
of $f^{\ot n}(TV\,\#\GG)f^{\ot n}$ by the two sided ideal generated
by $f^{\ot n} \C[\G^n]J \C[\G^n]f^{\ot n}$.
We will show that, via the isomorphism $T_\B \E \#S_n \simeq 
f^{\ot n}(TV\#\,\GG)f^{\ot n}$ constructed above, the two sided ideal 
that  $f^{\ot n} \C[\G^n]J \C[\G^n]f^{\ot n}$ generates 
gives precisely the relations that define $\A_{n,\la,\nu}$ as a
quotient of $T_\B \E \#S_n$.
We will use the description of $J$ given by Lemma \ref{leh}.

First, we consider the relations (R1) in Lemma \ref{leh}.
Observe that for any $g \in \G$, since $\om_L$ is $\G$-invariant,
we have
$$ g(x\ot y -y\ot x) = (x\ot y -y\ot x)g \,\in TL\,\#\G\,, $$ 
and since $c\in\ZZ \G$, we have
$$ g(t\cdot 1 + \sum_{\g \neq 1}c'_{\g}\g) =
 (t\cdot 1 + \sum_{\g \neq 1}c'_{\g}\g)g \,\in \C\G\,.$$
Also, for any $g,h \in \G$, 
$$ g_i h_j \big(\sum_{\g \in \G} s_{ij}\g_i\g_j^{-1}\big)
= \sum_{\g \in \G} s_{ij}(h\g)_i(\g g^{-1})_j^{-1}
= \big(\sum_{\g \in \G} s_{ij}\g_i\g_j^{-1}\big) g_i h_j \,\in 
\C[\GG]\,.$$
Hence, if $J_1 \sset J$ is spanned by elements of 
type (R1), then $$f^{\ot n}\C[\G^n]J_1\C[\G^n]f^{\ot n}
= f^{\ot n}J_1f^{\ot n}\C[\G^n]f^{\ot n}\,.$$
For any $i_1, \ldots, i_n\in I$ and $\ell \in [1,n]$, we have
\begin{align}
f_{i_{1}} \otot f_{i_{n}}\cdot [x_\ell, y_\ell] 
& = [x_\ell, y_\ell]\cdot f_{i_{1}} \otot f_{i_{n}} \nonumber\\
 = f_{i_{1}} \otot 
\frac{1}{\delta_{i_\ell}} & (\sum_{a\in Q, h(a)=i_\ell}
\phi_{a}\theta_{a} - \sum_{a\in Q, t(a)=i_\ell} \theta_{a}\phi_{a})
\otot f_{i_{n}} \nonumber
\end{align}
and 
$$f_{i_{1}} \otot f_{i_{n}}(t\cdot 1 + \sum_{\g \neq
1}c'_{\g}\g_\ell) = \frac{1}{\delta_{i_\ell}} \lambda_{i_\ell} f_{i_{1}}
\otot f_{i_{n}}\,.$$
Moreover, keeping in mind that we have fixed an 
isomorphism $\C\G\simeq \bigoplus_{i\in I}{\mathrm{Mat}}(\delta_i\times 
\delta_i)$, 
we have by orthogonality relations for matrix coefficients 
(see \cite[p.14]{Se}) 
\[ f_{i_{1}} \otot f_{i_{n}} \cdot
(\sum_{\g}s_{\ell j}\g_\ell\g^{-1}_{j})\cdot f_{i_{1}} \otot f_{i_{n}}
= \left\{ \begin{array}{ll}
\frac{|\G|}{\delta_{i_\ell}} f_{i_{1}}\otot f_{i_{n}} &
\textrm{if $i_j=i_\ell$}\,, \\
0 & \textrm{else}\,. \end{array} \right. \]
(Note that $f_i \g f_i$ (where $i\in I$,$\g\in \G$) is equal to
$f_i$ times the corresponding matrix coefficient for the
action of $\g$ on $N_i$.) 
Hence, (R1) gives the relations \vi of Definition \ref{dea}.

Next, we find the relations that come from (R2) in Lemma \ref{leh}.
To ease notations, we will assume without loss
of generality that $n=2$. (See Example \ref{ex2}.)

For any $u,v \in L$ and $g,h \in \G$, note that
$$ (g\ot h) \cdot [u_1,v_2] = [(gu)_1,(hv)_2]\cdot (g\ot h) \,,$$
and
\begin{align} 
(g\ot h) \cdot (\sum_{\g} \omega_{L}(\g u, v) s_{12}
\g_{1} \g_{2}^{-1}) 
&= \sum_{\g} \omega_{L}(\g u, v) s_{12}
(h\g)_{1} (\g g^{-1})_{2}^{-1} \nonumber\\
&= (\sum_{\g} \omega_{L}(\g gu, hv) s_{12}\g_{1}\g_{2}^{-1}) 
\cdot (g\ot h) \,.  \nonumber
\end{align}

Now, for any $i,j,k,l\in I$, we have
\begin{align} \label{eqn20}
(f_{i} & \otimes f_{j})(g\otimes 1) [u_{1}, v_{2}] (1 \otimes h)
(f_{k} \otimes f_{l})   \nonumber\\
=\,\, & \sum_{i_{1}, p_{1}, i_{2}, p_{2}}
(f_{i}g \otimes f_{j})\big(u_{1} \otimes (E^{i_{1}}_{p_{1},1} \otimes
E^{i_{2}}_{p_{2},1})\big) \bigotimes (E^{i_{1}}_{1, p_{1}} \otimes
E^{i_{2}}_{1, p_{2}})\big(v_{2} \otimes (f_{k} 
\otimes hf_{l})\big) \nonumber\\
-& \sum_{i_{1}, p_{1}, i_{2}, p_{2}}
(f_{i}g \otimes f_{j})
\big(v_{2} \otimes (g^{-1}E^{i_{1}}_{p_{1},1} \otimes
hE^{i_{2}}_{p_{2},1})\big) \bigotimes (E^{i_{1}}_{1, p_{1}}g \otimes
E^{i_{2}}_{1, p_{2}}h^{-1})\big(u_{1} \otimes 
(f_{k} \otimes hf_{l})\big) \nonumber\\
=\,\, & (f_{i}g\otimes f_{j})\big(u_{1} \otimes 
(f_{k} \otimes f_{j})\big)
\bigotimes (f_{k} \otimes f_{j})\big(v_{2} \otimes 
(f_{k} \otimes hf_{l})\big) \nonumber\\
 & - (f_{i}\otimes f_{j})\big(v_{2} \otimes 
(f_{i} \otimes hf_{l})\big)
\bigotimes (f_{i}g \otimes f_{l})\big(u_{1} 
\otimes (f_{k} \otimes f_{l})\big) 
\end{align}
Note that via our identifications in \S \ref{ss34}, 
the last line of (\ref{eqn20}) is an element of
$T^2_\B\E = \E \bigotimes_\B \E$.
Now, on the other hand,
\begin{align} \label{eqn21}
(f_{i} & \otimes f_{j}) (g\otimes 1)
\big(\sum_{\g} \omega_{L}(\g u, v) s_{12}
\g_{1} \g_{2}^{-1}\big) (1 \otimes h) (f_{k} \otimes f_{l}) \nonumber\\
= & s_{12}\sum_{\g} \omega_{L}(\g u, v) (f_{j}\g f_{k}) \otimes  
(f_{i} g\g^{-1}h f_{l}) \,.
\end{align}
Observe that for any edge $a \in \QQ$, we can find
$g_{a}, h_{a} \in \G$ and $u_{a}, v_{a} \in L$ such that
$$f_{t(a)}g_{a}(u_{a} \otimes f_{h(a)}) \neq 0 \qquad\textrm{and}\qquad
f_{h(a)}(v_{a} \otimes h_{a}f_{t(a)}) \neq 0 \,.$$
Suppose that $Q$ is not of type $A_1$, so that each of the spaces
$f_i(L \ot \C\G)f_j$ is at most one dimensional.
For any $i,j\in I$, there is a bijection
$$(f_{i}\mathbb{C}\G \otimes L\otimes \mathbb{C}\G f_{j})^{\G}
\rightarrow f_{i}(L\otimes \mathbb{C}\G)f_{j}:
\alpha\otimes u \otimes \beta \mapsto \alpha(u \otimes \beta)\,.$$
Here, the action of $\g\in\G$ on 
$\alpha\otimes u \otimes \beta\in f_{i}\mathbb{C}\G \otimes L\otimes 
\mathbb{C}\G f_{j}$ is $\alpha\g^{-1}\otimes \g u \otimes \g\beta$.
There is also a $\G$-equivariant non-degenerate pairing
$$(f_{i}\mathbb{C}\G \otimes L\otimes \mathbb{C}\G
f_{j}) \bigotimes (f_{j}\mathbb{C}\G \otimes L\otimes \mathbb{C}\G
f_{i}) \rightarrow \mathbb{C} $$
$$ (\alpha\otimes u \otimes \beta) \bigotimes
(\alpha'\otimes u' \otimes \beta') \mapsto
(\alpha\beta')(\alpha'\beta)\omega_{L}(u, u') $$
Thus, for any edge $a \in \QQ$, we may assume that 
$\omega_{L}(u_{a}, v_{a})=1$. Moreover, 
$f_{t(a)}g(v_{a} \otimes f_{h(a)}) = 0$
if $f_{h(a)}(v_{a} \otimes hf_{t(a)}) \neq 0$; and
$f_{h(a)}(u_{a} \otimes hf_{t(a)}) = 0$ if
$f_{t(a)}g(u_{a} \otimes f_{h(a)}) \neq 0$.

Note that if $k \neq j$ or $l \neq i$, then the expression in
(\ref{eqn21}) is zero. Hence, if $a: k \to i$ 
and $b: l \to j$ are two edges of $\QQ$ such
that $b \neq a^{*}$ or $a \neq b^*$, then we obtain from
(\ref{eqn20}) and (R2) that 
$$(a \otimes h(b))(t(a) \otimes b)- (h(a) \otimes b)(a \otimes t(b)) 
= 0 \,.$$

Now suppose $k=j$ and $l=i$. If $N_i$ is not an irreducible
component of $L\ot N_j$, then the expression in
(\ref{eqn21}) is zero by orthogonality of matrix coefficients
on $N_i$ and $L\ot N_j$; in this case, note that (\ref{eqn20})
is also zero. Thus, suppose there is an edge $a: i \to j$ in $\QQ$. 
We consider the case that $a\in Q$; the case $a\notin Q$ is 
completely similar.
We have $N_{i} \subset L \otimes N_{j}$ via $\theta_{a}$.
Hence, we have a decomposition into irreducible components
$L \otimes N_{j} = N_{i} \bigoplus \cdots$. 
A basis for $N_{i} \bigoplus \cdots$ compatible with 
this direct sum decomposition is
$$\xi_{1} := f_{i}= E^{i}_{1,1}\,,\quad \xi_{2} := E^{i}_{2,1}\,, 
\quad \xi_3 :=E^i_{3,1}\,, \quad\ldots \,\,,$$
and a basis for $L \otimes N_{j}$ is
$$\mu_{1} := u_a\otimes f_{j} = u_a\otimes E^{j}_{1,1}\,,
\quad \mu_{2} := u_a\otimes E^{j}_{2,1}\,, 
\quad\ldots\,\,,\quad \mu_{2 \delta_{j}}
:= v_a\otimes E^{j}_{\delta_{j}, 1} \,.$$
Define the matrix $\tau = (\tau_{p,q})$ by
$\mu_{q}= \sum_{p}\tau_{p,q}g^{-1}\xi_{p}$ and the matrix
$\varrho = (\varrho_{p,q})$ by $h\xi_q = 
\sum_p \varrho_{p,q} \mu_p$.

Using the fact that the composition
$$N_{j} \stackrel{\phi_{a}}{\hookrightarrow} L\otimes N_{i}
\stackrel{\theta_{a}}{\hookrightarrow} L\otimes L
\otimes N_{j} \stackrel{\omega_{L}\otimes 1}{\rightarrow} N_{j}$$
is multiplication by $\delta_{i}$, we get
$$f_{i}g(u_a \otimes f_{j})= \tau_{1,1}f_{i} \qquad \textrm{and}\qquad
f_{j}(v_a \otimes hf_{i}) 
= -\frac{\varrho_{1,1}}{\delta_{i}}f_{j}\,.$$

Now consider the matrix coefficients 
for $\g: L\ot N_j \map  L\ot N_j$,
where the matrix representing $\g$ is taken with respect to
the basis $\{g^{-1}\xi_p \}_{p=1,2,\ldots}$
for the domain of $\g$ and
the basis $\{h\xi_p \}_{p=1,2,\ldots}$
for the image of $\g$, vice versa
for the matrix representing $\g^{-1}$.
Since each irreducible representation of $\G$
appears in $ L\ot N_j$ at most once,
and the bases $\{g^{-1}\xi_p \}_{p=1,2,\ldots}$ \
and $\{h\xi_p \}_{p=1,2,\ldots}$ respect
the decomposition into irreducible components, the usual
orthogonality relations, proved for example in \cite[p.14]{Se},
continue to hold in this situation, that is, we have
$$ \sum_{\g \in \G}   
\omega_{L}(\g u_a,v_a) (f_{j}\g f_{j})\ot(f_{i}g \g^{-1} h f_{i})
= \frac{|\G| \tau_{1,1} \varrho_{1,1}}{\delta_{i}} f_j\ot f_i\,.$$
Hence, taking $u=u_a$ and $v=v_a$ in (\ref{eqn20})--(\ref{eqn21}),
we get from (R2) the relation:
$$ (a^{*}\otimes h(a))(h(a) \otimes a)
-(t(a)\otimes a)(a^{*} \otimes t(a))
= s_{12} \frac{k |\G|}{2}(h(a)\otimes t(a)) \,.$$
Note that $\tau_{1,1}$ and $\varrho_{1,1}$ are non-zero when
$g=g_a$ and $h=h_a$.

Now take $u=u_a$ and $v=u_a$ in (\ref{eqn20})--(\ref{eqn21}).
If $f_{i}g(u_a \otimes f_{j}) \neq 0$, then  
$f_{j}(u_a \otimes hf_{i})=0$ and so
$\varrho_{\delta_{j}+1,1}= 0$. 
Thus, both sides of the relation that (R2) gives in this case
are zero. Similarly if $u=v_a$ and $v=v_a$.

When $Q$ is of type $A_1$, $\G=\Z/2\Z$, and in particular
it is abelian. It is straightforward to check the relations 
directly in this case.

We conclude that (R2) gives the relations \vii of Definition
\ref{dea}.
\end{proof}

\setcounter{equation}{0}
\footnotesize{

Department of Mathematics, Massachusetts Institute of Technology,
Cambridge, MA 02139, USA;\\
\hphantom{x}\quad\, {\tt wlgan@math.mit.edu}

\smallskip

Department of Mathematics, University of Chicago, 
Chicago, IL 60637, USA;\\ 
\hphantom{x}\quad\, {\tt ginzburg@math.uchicago.edu}}


\begin{thebibliography}{APK}

\bibitem[BG]{BG} A. Braverman, D. Gaitsgory,
{\em Poincar\'e-Birkhoff-Witt theorem for quadratic algebras 
of Koszul type}, J. Algebra {\bf 181} (1996), no. 2, 315--328,
{\tt hep-th/9411113}.

\bibitem[CB]{CB} W. Crawley-Boevey,
{\em DMV Lectures on Representations of quivers, preprojective 
algebras, and deformations of quotient singularities}, 1999,
available at:
{\tt http://www.amsta.leeds.ac.uk/$\sim$pmtwc/}.

\bibitem[CBH]{CBH} W. Crawley-Boevey, M. Holland,
{\em Noncommutative deformations of Kleinian singularities},
Duke Math. J. {\bf 92} (1998), no. 3, 605--635. 

\bibitem[EG]{EG} P. Etingof, V. Ginzburg,
{\em Symplectic reflection algebras, Calogero-Moser space, and 
deformed Harish-Chandra homomorphism}, Invent. Math. {\bf 147} 
(2002), no. 2, 243--348, {\tt math.AG/0011114}.

\bibitem[EO]{EO} P. Etingof, A. Oblomkov,
{\em Quantization, orbifold cohomology, and Cherednik algebras},
preprint, {\tt math.QA/0311005}.

\bibitem[Gr]{Gr} E.L. Green,
{\em Introduction to Koszul algebras}, Representation theory 
and algebraic geometry (Waltham, MA, 1995), 45--62, 
London Math. Soc. Lecture Note Ser., 238, 
Cambridge Univ. Press, Cambridge, 1997.

\bibitem[GS]{GS} I. Gordon, S.P. Smith, 
{\em Representations of symplectic reflection algebras and 
resolutions of deformations of symplectic quotient singularities},
Math. Ann. {\bf 330}  (2004),  no. 1, 185--200, 
{\tt math.RT/0310187}.

\bibitem[MV]{MV} R. Martinez-Villa,
{\em Applications of Koszul algebras: the preprojective algebra},
Representation theory of algebras (Cocoyoc, 1994), 487--504, 
CMS Conf. Proc., 18, 
Amer. Math. Soc., Providence, RI, 1996. 

\bibitem[Se]{Se} J.-P. Serre,
{\em Linear representations of finite groups},
Graduate Texts in Mathematics, Vol. 42. Springer-Verlag, New 
York-Heidelberg, 1977.

\end{thebibliography}
\end{document}